\documentclass{cmslatex}
\usepackage[paperwidth=7in, paperheight=10in, margin=.875in]{geometry}
\usepackage[colorlinks,linkcolor=black,anchorcolor=black,citecolor=black,urlcolor=black]{hyperref}
\usepackage{amsfonts,amssymb}
\usepackage{amsmath}
\usepackage{graphicx}
\usepackage{cite}
\usepackage{enumerate}
\renewcommand{\theequation}{\arabic{section}.\arabic{equation}}

\usepackage{cleveref}
\usepackage{xspace}
\usepackage{subcaption}
\usepackage{enumitem}

\crefname{equation}{}{}

\newcommand{\B}{\ensuremath{\mathbf{B}}\xspace}
\newcommand{\C}{\ensuremath{\mathbf{C}}\xspace}
\renewcommand{\b}{\ensuremath{\mathbf{b}}\xspace}
\renewcommand{\c}{\ensuremath{\mathbf{c}}\xspace}

\newcommand{\M}{\ensuremath{\mathbf{M}}\xspace}
\newcommand{\m}{\ensuremath{\mathbf{m}}\xspace}

\renewcommand{\H}{\ensuremath{\mathbf{H}}\xspace}
\newcommand{\h}{\ensuremath{\mathbf{h}}\xspace}

\newcommand{\Real}{\ensuremath{\mathbb{R}}\xspace}

\newcommand{\F}{\ensuremath{\mathbf{F}}\xspace}
\newcommand{\f}{\ensuremath{\mathbf{f}}\xspace}
\newcommand{\A}{\ensuremath{\mathbf{A}}\xspace}
\newcommand{\I}{\ensuremath{\mathbf{I}}\xspace}

\renewcommand{\k}{\ensuremath{\mathbf{k}}\xspace}

\newcommand{\w}{\ensuremath{\mathbf{w}}\xspace}
\renewcommand{\v}{\ensuremath{\mathbf{v}}\xspace}
\renewcommand{\u}{\ensuremath{\mathbf{u}}\xspace}
\newcommand{\N}{\ensuremath{\mathbf{N}}\xspace}

\newcommand{\bnabla}{\boldsymbol\nabla}

\begin{document}

\title{A finite element based heterogeneous multiscale method for
  the Landau-Lifshitz equation} \author{Lena Leitenmaier
  \thanks{Department of Mathematics, KTH, Royal Institute of
    Technology, Stockholm, Sweden, (\texttt{lenalei@kth.se}).}  \and Murtazo
  Nazarov \thanks{Division of Scientific Computing, Department of
    Information Technology, Uppsala University, Sweden
    (\texttt{murtazo.nazarov@it.uu.se}).}
}

 \maketitle

\begin{abstract}
  We present a Heterogeneous Multiscale Method for the
  Landau-Lifshitz equation with a highly oscillatory diffusion
  coefficient, a simple model for a ferromagnetic composite. A
  finite element macro scheme is combined with a finite difference
  micro model to approximate the effective equation corresponding to
  the original problem. This makes it possible to obtain effective
  solutions to problems with rapid material variations on a small
  scale, described by $\varepsilon \ll 1$, which would be too
  expensive to resolve in a conventional simulation.
\end{abstract}



\section{Introduction}
Micromagnetic simulations of ferromagnetic materials provide an
important tool in physics and material science. The dynamics of the
magnetization $\M^\varepsilon: \Omega \times (0, T] \to \Real^3$ are typically described using the
Landau-Lifshitz equation,
\begin{subequations} \label{eq:LL}
\begin{align}
  \partial_t \M^\varepsilon(x,t) &= - \M^\varepsilon \times \H^\varepsilon(\M^\varepsilon) - \alpha \M^\varepsilon \times (\M^\varepsilon \times \H^\varepsilon(\M^\varepsilon)),  && x \in \Omega, \,t > 0,\\
  \M^\varepsilon(x, 0) &= \M_\mathrm{init}(x), && x \in \Omega, \,t = 0\\
  \bnabla \M^\varepsilon \cdot \mathbf{n} &= 0,  &&\,\,x \in \partial \Omega,\, t > 0,
\end{align}
\end{subequations}
where $\alpha$ is a material dependent parameter determining the
strength of damping and the initial data $\M_\mathrm{init}$ is such
that $|\M_\mathrm{init}| = 1$ throughout $\Omega$. The vector
$\mathbf{n}$ is the normal to the boundary $\partial \Omega$.
Moreover, $\H^\varepsilon$ denotes the effective field affecting the
magnetization.  In this paper, we consider the case of a
ferromagnetic composite. A simplified model for this is to introduce
a material coefficient $a^\varepsilon$ describing the variations in
the material, which are on a scale characterized by the parameter
$\varepsilon \ll 1$. 
This type of description has been used in several approaches
recently, for example in
\cite{alouges2019stochastic,alouges2015homogenization,multilayer,highcontrast}.
With this model, the effective field $\H^\varepsilon$ we consider is
\[\H^\varepsilon(\M^\varepsilon) := \bnabla \cdot (a^\varepsilon \bnabla
  \M^\varepsilon) + \H_\mathrm{low}(\M^\varepsilon),\]
where the first term is due to the exchange interaction between
magnetic moments in the material, influenced by the material
coefficient, while $\H_\mathrm{low}$ represents lower order terms,
in particular external field, anisotropy and the so-called
demagnetization field,
\[\H_\mathrm{low}(\M^\varepsilon) = \H_\mathrm{ext} + \H_\mathrm{ani}(\M^\varepsilon) +   \H_\mathrm{dem}(\M^\varepsilon).\]
For small values of $\varepsilon$, direct numerical simulation of
\cref{eq:LL} is infeasible since the computational cost becomes too
high when resolving the $\varepsilon$-scale. We therefore use the
framework of Heterogeneous Multiscale Methods (HMM)
\cite{weinan2003,acta_numerica}, which makes it possible to
numerically obtain an approximation to the effective solution to the
problem. The idea with this framework is to combine a coarse scale
macro model, involving a missing quantity that encodes the effect of
the fast variations, with a micro model that resolves the fine
scale. The micro model is only solved on a small domain in time and
space, keeping the computational cost independent of the scale of
the fast variations. The solution to the micro model is then used to
approximate the unknown quantity required to complete the macro
model.

Several ways to set up HMM for a periodic version of \cref{eq:LL}
are discussed in \cite{paper2}. A finite difference based
implementation of one of these approaches, the so-called field
model, is studied in \cite{paper3}. In this article, we focus on a variation of the
so-called flux model and investigate how to combine a finite element
macro scheme with a finite difference discretized micro model.  This
makes it possible to use the approach for more general geometries
and gives a high flexibility. Additionally, in contrast to
\cite{paper3}, the effective field $\H^\varepsilon$ considered in this paper is
more general and contains not only the exchange term but also
applied field and demagnetization.

This article is structured as follows. We first introduce useful
notation and give a definition of the finite element spaces used in
\Cref{sec:prelim}. In \Cref{sec:HMM}, the HMM approach in general as
well as the considered HMM macro and micro model and their numerical
solution are described. Related error estimates are given and
illustrated with an example. Finally, in \Cref{sec:num_ex},
numerical examples are given to demonstrate the properties of the
scheme.

\section{Preliminaries}\label{sec:prelim}

Throughout this article, we consider a domain
$\Omega \subset \Real^d$, where $d = 2$ or $3$. For numerical
examples, we use $d = 2$.
We let $\bnabla \m$ denote the Jacobian of $\m \in \Real^3$,
\[\bnabla \m = [\partial_{x_1} \m \cdots \partial_{x_d} \m].\]
Furthermore, we use the colon-operator to denote the column-wise
scalar product of two matrices. Consider matrices
$\B, \C \in \Real^{3 \times d}$ with columns
$\b_j, \c_j \in \Real^3$, respectively, then
\[\B : \C = \sum_{j=1}^d \b_j \cdot \c_j.\]
In general, we use the convention that scalar and cross product
between a vector and a matrix are done column-wise, and that scalar
differential operators are applied element-wise to vector-valued
functions. Furthermore, the divergence operator is applied row-wise
to a matrix-valued function in order to have consistency with the
scalar case in the sense that
\[\bnabla \cdot (\bnabla \m) =  \Delta \m.\]

For the finite element discretization, we introduce an affine mesh
$\mathcal{T}_h$ which is a subdivision of $\Omega$ into disjoint
elements $K$ such that
\[\overline{\Omega} = \bigcup_{K \in \mathcal{T}_h} \overline{K},\]
where $\overline{\Omega}$ and $\overline{K}$ denote the closures of
$\Omega$ and $K$, respectively. We consider a family of
shape-regular meshes, $\{\mathcal{T}_h\}_{h>0}$, such that each mesh is
conforming. The shortest edge in a given triangulation is denoted
$H_\mathrm{min}$.

Throughout this paper, we consider Lagrange finite elements and
denote the set of nodes $\{N_1, ..., N_J\} =: \mathcal{N}_h$. The
associated pieceswise linear scalar nodal basis functions are
$\{\phi_j(x)\}_{j = 1}^J$, defined such that
$\phi_j(x_i) = \delta_{ji}$ for any $i, j \le J$. Moreover, let
$\chi_K$ be indicator functions for $K \in \mathcal{T}_h$. Then
the space of piecewise linear vector-valued functions is given by
\begin{equation}
  \label{eq:vh_space}
  V_h := \{\v_h \in \mathcal{C}^0(\overline \Omega; \Real^3) \, |\, \v_h(x) = \sum_{j = 1}^J \v_j
  \phi_j(x), ~\text{where}~ \v_j \in \Real^3\},
\end{equation}
and the space of piecewise constant vector-valued functions is
\begin{equation}
  \label{eq:wh_space}
  W_h := \{\w_h \,|\, \w_h(x) = \sum_{K \in \mathcal{T}_h} \w_K \chi_K,  ~\text{where}~ \w_K \in \Real^3\}.
\end{equation}
We define the interpolation operator $\mathcal{I}_h : \mathcal{C}^0(\Omega; \Real^3) \to V_h$ such that
\[\mathcal{I}_h(\m) := \sum_{j=1}^J \m_j \phi_j(x), \quad \text{where} \quad\m_j := \m(N_j).\]
Note that the Landau-Lifshitz equation \cref{eq:LL} is length
preserving, due to its cross product structure it holds that
\begin{align}\label{eq:length_pres}
  \partial_t |\M^\varepsilon|^2 = 2 \M^\varepsilon \cdot \partial \M^\varepsilon = 0.
\end{align}
Hence $|\M_\mathrm{init}| = 1$ implies that
$|\M^\varepsilon(x, t)| = 1$ for all $x \in \Omega$ and
$0 \le t \le T$. To accommodate this normalization constraint in the
finite element solution, we introduce the solution space
\begin{equation}
  \label{eq:mh_space}
  M_h := \{\m_h \in V_h \,|\, \m_h(x) = \sum_{j = 1}^J \m_j
  \phi_j(x) ~\text{with}~|\m_j| = 1 \}.
\end{equation}

To make it easier to distinguish between solutions to the HMM micro
and macro model, we in general use capitals, for example
$\M, \M^\varepsilon$, to refer to solutions on the whole domain
$\Omega$. To denote solutions on micro domains, we use lowercase
letters, such as $\m^\varepsilon$.

\section{Heterogeneous Multiscale Methods}\label{sec:HMM}

The framework of Heterogeneous Multiscale Methods (HMM) was first
introduced by Engquist and E in \cite{weinan2003}. The goal with the
approach is to achieve numerical homogenization for multiscale
problems with scale separation. To accomplish this, one combines a
macro and micro model in such a way that the relevant influence of
the fast variations in the problem is captured in the micro problem
and encoded by an effective quantity, that can then be used to solve
the macro scheme on a rather coarse discretization. For a wide range
of applications, HMM has been shown to be an efficient way to obtain
effective solutions to multiscale problems, as for example described
in \cite{acta_numerica,hmm1}. Note that due to the fact that an
effective solution is approximated by HMM, some error compared to a
(numerical) solution to the original problem is introduced. The size
of this homogenization error is typically determined by the scale of
the fast variations, which are not included in the HMM macro
solution.  In \cite{paper1}, it was proved that for strong solutions
to \cref{eq:LL} in a periodic setting and corresponding effective
solutions, this error is $\mathcal{O}(\varepsilon)$.

\subsection{HMM for the Landau-Lifshitz equation}

To set up a HMM scheme for the Landau-Lifshitz problem \cref{eq:LL},
consider first the case of a periodic material coefficient,
$a^\varepsilon(x) = a(x/\varepsilon)$, where we assume that
$a(y) \in \mathcal{C}^\infty(\Omega)$ is bounded by positive constants
$a_\mathrm{min}$ and $a_\mathrm{max}$,
$0 < a_\mathrm{min} \le a(y) \le a_\mathrm{max}$ for all $y \in
\Omega$.
As we moreover have that the initial data for \cref{eq:LL},
$\M_\mathrm{init}$, is such that $|\M_\mathrm{init}(x)| = 1$ for all
$x \in \Omega$, independent of the material coefficient, 
we conclude based on
\cite{paper1} that the homogenized problem corresponding to
\cref{eq:LL} is to find $\M_0: \Omega \times [0, T] \to \Real^3$ such
that for $x \in \Omega$ and $0 \le t \le T$,
\begin{subequations}\label{eq:hom}
\begin{align}
  \partial_t \M_0(x,t) &= - \M_0 \times \left[\bnabla \cdot (\bnabla \M_0 \A^H) + \H_\mathrm{low}^\mathrm{hom}(\M_0)\right] \label{eq:hom_a}
                  \\ &\hspace{1cm}- \alpha \M_0 \times \left[\M_0 \times (\bnabla \cdot (\bnabla \M_0 \A^H) + \H_\mathrm{low}^\mathrm{hom}(\M_0))\right] , \nonumber \\
  \M_0(x,0) &= \M_\mathrm{init}(x), \\
  (\bnabla \M_0 \A^H) \cdot \mathbf{n} &= 0 \quad \text{on } \partial \Omega.
\end{align}
\end{subequations}
The initial data for \cref{eq:hom} is the same as for \cref{eq:LL}
and $\H_\mathrm{low}^\mathrm{hom}$ denotes a homogenized version of
the lower order field terms. The homogenized coefficient matrix
$\A^H$ for the periodic case is the same as for standard elliptic
homogenization problems and can be computed by solving
\begin{align}\label{eq:AH}
  \A^H = \int_{[0, 1]^d} a(y) (\I + (\bnabla_y \boldsymbol \chi(y))^T) dy,
\end{align}
where $\I$ denotes the $d\times d$ identity matrix and
$\boldsymbol \chi$ solves the so-called cell problem
\[\bnabla \cdot(a(y) \bnabla \boldsymbol \chi(y)) = - \nabla_y a(y).\]

\begin{remark}\label{remark1}
  In the rest of this article, we do not consider problems including
  anisotropy effects. We moreover assume that the variations in the
  material primarily affect the exchange term, since this is a very
  short range interaction. The demagnetization, in contrast, is a
  long range effect. The applied field is in general independent of
  the magnetization itself. We hence choose in the following to
  use the approximation
  \[\H_\mathrm{low}^\mathrm{hom}(\M) \approx \H_L(\M) := \H_\mathrm{ext} +
    \H_\mathrm{dem}(\M),\] which is included in the HMM macro model
  only.  Neglecting changes in the so-called saturation
  magnetization between the materials, this approximation matches
  with the results in
  \cite{alouges2019stochastic,alouges2015homogenization}.
\end{remark}

Since the homogenized matrix $\A^H$ is symmetric, we have the identity
\[\M_0 \times \bnabla \cdot (  \bnabla \M_0  \A^H) = \bnabla \cdot (\M_0 \times (\bnabla \M_0 \A^H))\,,\]
and by the
vector triple product identity together with the fact that $|\M_0| = 1$,
it follows that
\[- \M_0 \times \left[\M_0 \times \bnabla \cdot ( \bnabla \M_0 \A^H)\right] = \bnabla
  \cdot (\bnabla \M_0 \A^H) + (\bnabla \M_0 : (\bnabla \M_0 \A^H)) \M_0.\]
Hence, \cref{eq:hom_a}
can be rewritten as
\begin{align} \label{eq:hom_b}
  \partial_t \M_0 &= - \bnabla \cdot (\M_0 \times (\bnabla \M_0 \A^H)) + \alpha \left[\bnabla
                  \cdot (\bnabla \M_0 \A^H) + (\bnabla \M_0 : (\bnabla \M_0 \A^H)) \M_0\right]
                  \\ &\quad - \M_0 \times (\H_L(\M_0) + \alpha \M_0 \times \H_L(\M_0)). \nonumber
\end{align}
Taking \cref{eq:hom_b} as an inspiration, we deduce that a possible
HMM macro model for the problem \cref{eq:LL} is to find $\M(x,t)$ such
that for $x \in \Omega$ and $0 \le t \le T$,
\begin{subequations}\label{eq:macro_strong}
\begin{align}
  \partial_t \M &= - \bnabla \cdot(\M \times  \F) + \alpha \left[\bnabla \cdot \F + (\bnabla \M : \F) \M\right]
                    - \M \times \left[\H_L(\M)  + \alpha \M \times \H_L(\M)\right],  \\
  \M(x,0) &= \M_\mathrm{init}(x), \\
  \F \cdot \mathbf{n} &= 0 \quad \text{on } \partial \Omega, \label{eq:neumann_bc}
\end{align}
\end{subequations}
where the flux $\F(\M, x)$ is unknown and has to be approximated at
each discrete point in time $t^k$ where it is needed in a numerical
scheme. In case of a periodic material coefficient,
$\F(\M, x) \approx \bnabla \M \A^H$.

The approximation of $\F$ at a point $x \in \Omega$ is
based on the solution $\m^\varepsilon$ to the micro problem
\begin{subequations}\label{eq:micro}
\begin{align}
  \partial_t \m^\varepsilon(\xi,\tau) &= - \m^\varepsilon \times \bnabla \cdot (a^\varepsilon \bnabla \m^\varepsilon)
  - \alpha \m^\varepsilon \times \left[\m^\varepsilon \times \bnabla \cdot (a^\varepsilon \bnabla \m^\varepsilon)\right], \\
  \m^\varepsilon(\xi, 0) &= \m_\mathrm{init}(\xi) := \M^k(x + \xi),
\end{align}
\end{subequations}
for $\xi \in \Omega_\mathrm{mic}$ and $0 \le \tau \le \eta$, where
$\eta \sim \varepsilon^2$ and the micro domain is
$\Omega_\mathrm{mic} = [-\mu', \mu']^d$ with
$\mu' \sim \varepsilon$. This implies that the size of the micro domain and time
interval is chosen proportional to the scale of the fast variations
in the problem \cite{paper1}.  The macro and micro problem are
coupled via the initial data to the micro problem, which is set
according to the current macro solution $\M^k(x) \approx \M(x, t^k)$
at a given discrete point in time $t^k$.  The solution of the micro
problem is described in more detail in \Cref{sec:upscaling}.

\vspace{.1cm}
\begin{remark}
  Note that in the micro problem \cref{eq:micro}, only the exchange
  contribution to the effective field is considered. We choose this
  model based on the considerations given in \cref{remark1}, and due
  to the fact that the micro problem only is solved on a small,
  local domain and for a short time interval. Hence we suppose that
  strong short-range exchange forces dominate all other forces
  here. An alternative approach would be to include the lower order
  field terms in the micro model as well. However, this would result
  in an increased computational cost as in particular the
  computation of the long-range demagnetization term is rather
  computationally expensive \cite{abert2013numerical}.
\end{remark}

Once the micro problem is solved, the quantity
$a^\varepsilon \bnabla \m^\varepsilon$ is averaged in space and time
to approximate $\F(\M^k)$. To reduce the approximation error
introduced in this process, we use smooth averaging kernels $k$ from
the space of kernels $\mathbb{K}^{p, q}$, see \cite{stiff,
  doghonay1}. This space of smoothing kernels is defined such that
$k \in \mathbb{K}^{p, q}$ given that
\[\int_{-1}^1 k(x) x^r dx =
  \begin{cases}
    1\,, & r = 0\,,\\
    0\,, & 1 \le r \le p\,,
  \end{cases}\]
and additionally,
\[k \in C_c^{q}(\Real) ~\text{with} ~ supp(k) = [-1, 1], ~k^{(q+1)} \in BV(\Real).\]
Moreover, in \cite{paper2} the subspace $\mathbb{K}_0^{p, q} \subset \mathbb{K}^{p, q}$
is defined such
\[\mathbb{K}_0^{p, q} := \{k \in \mathbb{K}^{p, q} | k(x) = 0 \text{~for~} x \le 0\}.\]
Following the conventions in the field, we use the notation
that $k_\mu(x)$ is a scaled version of $k$,
\[k_\mu(x) := \frac{1}{\mu} k(x/\mu).\]
In several space dimensions, $d > 1$, we let
\[k(x) := k(x_1) \cdot \ldots \cdot k(x_d).\]
Let now $\m^\varepsilon(\xi, \tau)$ be the solution to \cref{eq:micro}
with initial data $\m_\mathrm{init}(\xi) = \M^k(x + \xi)$ for $\xi \in \Omega_\mathrm{mic}$
and $0 \le \tau \le \eta$ with $\eta \sim \varepsilon^2$.
Then we define
\begin{align}\label{eq:flux}
  \F(\M^k, x) := \int_{[-\mu, \mu]^d} \int_0^\eta k_\mu(\xi) k_\eta^0(\tau) a^\varepsilon(\xi) \bnabla \m^\varepsilon(\xi, \tau) d\tau d\xi,
\end{align}
where $k \in \mathbb{K}^{p_x, q_x}$ and $k^0 \in \mathbb{K}^{p_t, q_t}_0$,  and the averaging parameter $\mu$ is
chosen such that $\mu \sim \varepsilon$, $\mu \le \mu'$.

Note that \cref{eq:macro_strong,eq:micro,eq:flux} is a variation of the so-called
flux model in \cite{paper2}. There estimates for the error
introduced when approximating $\F$ are given for the case of a
periodic material coefficient and under certain regularity
conditions, as well as under the assumption that \cref{eq:micro} is
solved throughout $\Omega$ rather than only on
$\Omega_\mathrm{mic}$.  Then the approximation error is bounded as
given in the following theorem from \cite{paper2}.

\vspace{.1cm}
\begin{theorem}\label{thm:paper2}
  Assume $a(y) \in \mathcal{C}^\infty(\Omega)$ such that
  $a_\mathrm{min} \le a(y) \le a_\mathrm{max}$,
  $0 < \varepsilon < 1$ and $0 < \alpha \le 1$ and let
  $\varepsilon^2 < \eta \le \varepsilon^{3/2}$. Suppose that
  $\m^\varepsilon(\xi, \tau) \in \mathcal{C}^1([0, \eta]; H^{2}(\Omega))$ is the
  exact solution to the micro problem \cref{eq:micro} for
  $0 \le \tau \le \eta$ and $\xi \in \Omega_\mathrm{mic} = [0, 1]^d$ with
  periodic boundary conditions.  Moreover, suppose that there is a
  constant $c$ independent of $\varepsilon$ such that
  $\|\bnabla \m^\varepsilon(\cdot, \tau)\|_{L^\infty} \le c$ and that
  the solution to the corresponding homogenized problem is
  $\m_0 \in \mathcal{C}^\infty(0, \eta, H^\infty(\Omega))$.
  Consider averaging kernels   $k \in \mathbb{K}^{p_x, q_x}$ and
  $k^0 \in \mathbb{K}^{p_t, q_t}_0$ and let  $\varepsilon < \mu < 1$.
  Then $\F$ as given by \cref{eq:flux} satisfies
  \begin{align*}
    &\left| \F - \bnabla \m_\mathrm{init}(0) \A^H \right|
      =: E_\varepsilon + E_\mu + E_\eta,
  \end{align*}
  where
  \begin{align}\label{eq:err_terms}
    E_\varepsilon \le C \varepsilon, \quad E_\mu \le C \left(\mu^{p_x +1} +  \left(\frac{\varepsilon}{\mu}\right)^{q_x + 2}\right) \quad \text{and} \quad E_\eta \le C \left( \eta^{p_t + 1} + \left(\frac{\varepsilon^2}{ \eta}\right)^{q_t+1}\right).
  \end{align}
  The constant $C$ is independent of $\varepsilon$, $\mu$ and $\eta$
  but depends on $K$, $K^0$ and $\alpha$.
\end{theorem}



\subsection{Macro model}

The macro problem \cref{eq:macro_strong} is solved numerically using the
finite element method.  To obtain a weak formulation for the
problem, we multiply \cref{eq:macro_strong} by a test function
$\w \in H^1(\Omega, \Real^3)$ and integrate over
$\Omega$. This yields the weak problem
\begin{align}\label{eq:macro_weak}
  \int_{\Omega} \partial_t \M \cdot \w dx &= \int_{\Omega} (\M \times \F) : \bnabla \w dx - \alpha \int_{\Omega} \F: \bnabla \w dx + \alpha \int_{\Omega} (\bnabla \M : \F) (\M \cdot \w) dx \\
  &\qquad - \int_{\Omega} \left(\M \times [\H_L(\M) + \alpha \M \times \H_L(\M)]\right) \cdot \w dx,  \nonumber
\end{align}
for every test function $\w \in H^1(\Omega, \Real^3)$. The boundary
terms here vanish due to the given homogeneous Neumann boundary
condition, \cref{eq:neumann_bc}.

Note that finite element schemes based on several different weak
formulations have been proposed for the Landau-Lifshitz
equation. For instance, the schemes in \cite{bartels2006convergence}
and \cite{alouges2006fe,alouges2008new}, derived from the so-called
Gilbert form of \cref{eq:LL}, an equivalent way to rewrite the
equation, are commonly used and can be seen as advantageous with
regard to mimicking certain physical properties of the continuous
problem. However, in this article we choose the more direct approach
in \cref{eq:macro_weak} since our main goal is to study the
combination of micro and macro model and the influence of the flux
$\F$ introduced by the HMM approximation. It is possible to
introduce the unknown flux $\F$, approximated from the micro
problem, in other weak formulations as well in a similar way.

To solve \cref{eq:macro_weak} numerically, we adopt an idea from
\cite{alouges2006fe} and introduce an unknown variable
$\v \approx \partial_t \M$ which replaces $\partial_t \M$ in
\cref{eq:macro_weak}. We furthermore discretize in space and time
and denote the discretized magnetization at time $t^k$ by
$\M_h^k(x) \approx \M(x, t^k)$. For each discrete time step, the problem then
becomes to find $\v_h \in V_h$ such that for all $\w_h \in V_h$,
\begin{align}\label{eq:macro_weak_discrete}
  \int_{\Omega} \v_h \cdot \w_h dx &= \int_{\Omega} (\M^k_h \times \F_h^k) : \bnabla \w_h dx + \alpha \int_{\Omega} (\bnabla \M_h^k : \F_h^k) (\M_h^k \cdot \w_h) dx \\
  &\quad- \alpha \int_{\Omega} \F_h^k: \bnabla \w_h dx  - \int_{\Omega} \left(\M_h^k \times [\H_L(\M^k_h) + \alpha \M^k_h \times \H_L(\M_h^k)]\right) \cdot \w_h dx, \nonumber
\end{align}
where $\M_h^k \in M_h$ is given and $\F_h^k \in W_h$ is computed by
solving the micro problem \cref{eq:micro} and averaging according to \cref{eq:flux}.
Using the notation $\M_j^k$, $j = 1,..., J$ for the nodal values of $\M_h^k$
as in \cref{eq:mh_space}, the time update then is given by
\begin{subequations}\label{eq:ts_euler}
\begin{align}
  \tilde \M_j^{k+1} &= \M_j^k + \Delta t \v_j, \qquad j = 1, ..., J, \\
  \M^{k+1}_h &= \sum_{j=1}^J \M_j^{k+1} \phi_j, \qquad \text{where}\quad \M_j^{k+1} = \tilde \M_j^{k+1} / |\tilde \M_j^{k+1}|,
\end{align}
\end{subequations}
where $\phi_j(x)$ are the piecewise linear basis functions as in the definition of $M_h$,
\cref{eq:mh_space}.
Since \cref{eq:macro_weak_discrete} is a linear problem, it can also be
formulated in terms of a bilinear form $b$ and a linear form $L$. Let
\begin{align}
  b(\u, \w) :&= \int_\Omega \u \cdot \w dx, \\
  L(\w; \M_h^k, \F_h^k) :&= \int_{\Omega} (\M^k_h \times \F_h^k) : \bnabla \w dx + \alpha \int_{\Omega} (\bnabla \M_h^k : \F_h^k) (\M_h^k \cdot \w) dx \\
  &\quad- \alpha \int_{\Omega} \F_h^k: \bnabla \w dx  - \int_{\Omega} \left(\M_h^k \times [\H_L(\M^k_h) + \alpha \M^k_h \times \H_L(\M_h^k)]\right) \cdot \w dx \nonumber,
\end{align}
then \cref{eq:macro_weak_discrete} is to find $\v_h \in V_h$ such that for all $\w_h \in V_h$,
\[b(\v_h, \w_h) = L(\w_h; \M_h, \F_h).\]

\begin{figure}[h]
  \centering
  \includegraphics[width=.6\textwidth]{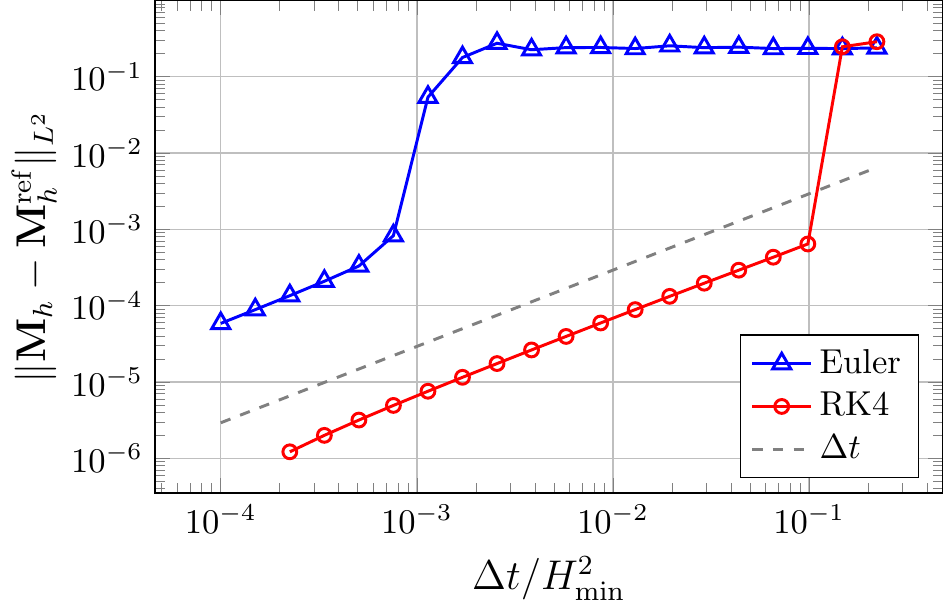}
  \caption{Approximation error with respect to a reference solution
    when solving $\Cref{eq:hom}$ for an example problem where
    $\alpha = 0.01$ with the finite element scheme described above
    for varying $\Delta t$. Reference solution $\M^\mathrm{ref}_h$
    computed with low $\Delta t$ on the same triangulation as
    $\M_h$. The considered example problem is the one described in
    \Cref{sec:ex1}, where more details are given. Here
    $H_\mathrm{min} = 0.171$ and the selected final time is $T = 1$.}
  \label{fig:stability}
\end{figure}

In the literature, for example
\cite{alouges2008new,bartels2006convergence}, it is shown that it in
many cases is advantageous to use implicit time integration rather
than an explicit Euler-like scheme as in \cref{eq:ts_euler}, since
with the latter, the time step size $\Delta t$ has to satisfy a
severe time step restriction to obtain stable
approximations. However, since the flux $\F$ is unknown and
approximated from the micro problem, which implies a very
complicated dependence on $\M$, it can in practice only be treated
explicitly in time. To obtain a slightly less harsh time step
restriction, we suggest to use a Runge-Kutta based time stepping
scheme instead of \cref{eq:ts_euler}. The overall scheme then is
described by the following steps:

\begin{itemize}
\item Initially, set
  \[\M_h^0 = \mathcal{I}_h(\M_\mathrm{init}).\]
\item Let $\Delta t = T/N$, where the number of time steps $N$ is
  such that $\Delta t < C \Delta H_\mathrm{min}^2$ for a constant
  $C < 1$, depending on $\alpha$. For $k = 0, ..., N-1$:
\begin{enumerate}
\item Let $\tilde \M = \M_h^k$
\item For $\ell = 1, .., 4$:
  \begin{itemize}
  \item for each face $K \in \mathcal{T}_h$, solve a micro problem
    \cref{eq:micro} with initial data given by $\tilde \M$ and compute the local
    value for the flux, $\F_K$, according to  \cref{eq:flux}.
    Then $\tilde \F_h \in W_h$ is
    \[\tilde \F_h(x) =  \sum_{K \in \mathcal{T}_h}  \F_K \chi_K,\]
    where $\chi_K$ are piecewise constant indicator functions as in the definition of $W_h$,
    \cref{eq:wh_space}.
  \item obtain $\k^\ell_h$ by solving
    \[b(\k^\ell_h, \w_h) = L(\w_h; \tilde \M, \tilde \F_h) \qquad \forall \w_h \in V_h.\]
  \item set $\tilde \M = \M_h^k + c_\ell \Delta t \k^\ell$,
    where $c_\ell = 1/2$ for $\ell = 1, 2$ and $c_\ell = 1$ for $\ell > 2$ according to the Butcher tableau.
  \end{itemize}
\item Use the stage value functions $\k_h^\ell$ to obtain $\tilde \M^{k+1}$,
  \[\tilde \M^{k+1} = \M^k_h + \frac{\Delta t}{6} \left(\k_h^1 + 2 \k_h^2 + 2 \k_h^3 + \k_h^4\right).\]
\item \label{item:renorm} Compute $\M^{k+1}_h$ using the normalized nodal values of $\tilde \M^{k+1}$,
  \[\M^{k+1}_h = \sum_{j=1}^J \M_j^{k+1} \phi_j, \qquad \text{where}\quad \M_j^{k+1} = \tilde \M_j^{k+1} / |\tilde \M_j^{k+1}|.\]
\end{enumerate}
\end{itemize}
This is what is applied for the numerical experiments in this
article. The effect of using a Runge-Kutta rather than Euler based
time step update on the time discretization error can for an example
problem be seen in \Cref{fig:stability}. Both stability and accuracy
of the approach improve when using a Runge-Kutta based time step, by
far making up for the additional computational effort due to the
stage value calculations. However, both approaches are only first
order accurate. This is due to the renormalization, step
\ref{item:renorm} in the description above, see for example also \cite{alouges2014precise}.
Note, though, that due to the given time step restriction,
$\Delta t < C H_\mathrm{min}^2$, we overall still have second order
accuracy with respect to the space discretization size.

\clearpage
\subsection{Micro problem solution and upscaling}\label{sec:upscaling}
\begin{figure}[h]
  \centering
  \includegraphics[width=.7\textwidth]{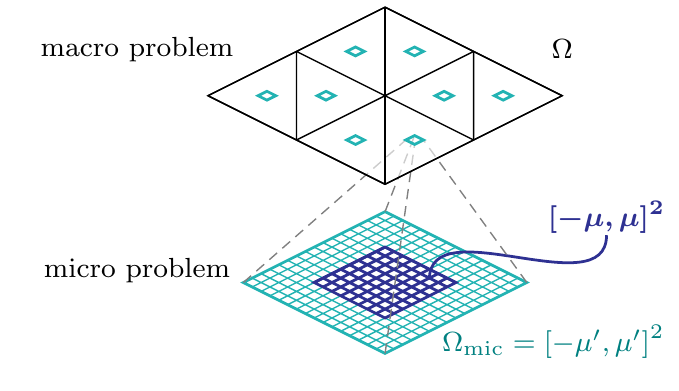}
  \caption{Domains involved in the HMM approach.}
  \label{fig:domains}
\end{figure}

~\\~
In the HMM micro problem around the macro location $x$ and at time
$t^k$, we aim to find $\m^\varepsilon(\xi, \tau)$ such that
\begin{subequations}\label{eq:micro2}
\begin{align}
  \partial_t \m^\varepsilon(\xi,\tau)
  &= - \m^\varepsilon \times \left[\bnabla \cdot (a^\varepsilon \bnabla \m^\varepsilon)
    + \alpha \m^\varepsilon \times \bnabla \cdot (a^\varepsilon \bnabla \m^\varepsilon)\right],
      &&\xi \in \Omega_\mathrm{mic}, \,0 < \tau \le \eta, \\
    \m^\varepsilon(\xi, 0) &= \m_\mathrm{init}(\xi) := \M^k(x + \xi), && \xi \in \Omega_\mathrm{mic}, ~\,\tau = 0, \\
  \m^\varepsilon(\xi, t) &= \m_\mathrm{init}(\xi),  && \xi \in \partial \Omega_\mathrm{mic}, \tau > 0,
\end{align}
\end{subequations}
where the considered final time is $\eta \sim \varepsilon^2$ and the
micro problem domain is $\Omega_\mathrm{mic} = [-\mu', \mu']^d$ with
$\mu' \sim \varepsilon$. We moreover suppose that $\mu'$ is such
that the whole micro domain $\Omega_\mathrm{mic}$ is inside one
triangle of the macro discretization, located around the
barycenter. Note that the latter assumption is due to reasons of
simplicity and comes naturally for small values of
$\varepsilon$. In the following, we focus on the case $d=2$. A schematic
overview of the connection between macro domain $\Omega$, micro
problem domain $\Omega_\mathrm{mic}$ and averaging domain
$[-\mu, \mu]^2$ is given in \Cref{fig:domains}.

The initial data $\m_\mathrm{init}$ for the micro problem is chosen
according to the current macro solution $\M^k_h$ at time $t^k$ in
the domain $\Omega_\mathrm{mic}$.  As $\M^k_h$ is a piecewise
bilinear function, the same holds for $\m_\mathrm{init}$.  Note that
this implies that $\m_\mathrm{init}$ is not normalized.  As a
consequence, an additional error is introduced in the upscaling
process.  In numerical experiments, this term appears to be
$\mathcal{O}(H)$, where $H$ describes the shortest edge of the
corresponding macro triangle.  In the following, we account for this
error by adding an additional term, $E_\mathrm{norm}$, to the error
estimate as stated in \Cref{thm:paper2}. Note that in the periodic
case and given a finite element approximation $\M_h \in M_h$ to an (unknown) actual solution $\M_\mathrm{ref}$, we
know that with the given scheme, it holds that
$\F(x, \M_h) \approx \bnabla \M_h(x) \A^H$, and that when using
piecewise linear finite elements,
$\bnabla \M_h = \bnabla \M_\mathrm{ref} + \mathcal{O}(H)$. Hence the
additional upscaling error is of the same order as this
approximation error.

Furthermore, as stated in \cref{eq:micro2}, the micro problem is
completed with Dirichlet boundary conditions. These boundary
conditions are artificial and cause an additional error, in the
following denoted by $E_\mathrm{\mu'}$. When considering an infinite
domain, that is as $\mu' - \mu \to \infty$, this error term
vanishes.  In practice, we cannot have an infinitely large
computational domain. We instead choose $\mu'$ such that
$E_\mathrm{\mu'}$ does not significantly influence the upscaling
error but not much larger in order to not increase the computational
cost more than necessary.  In general, $\mu'$ has to be chosen
larger when larger final times $\eta$ are considered, since with
increasing time, the errors caused by the artificial boundary
conditions travel further into the domain.

To solve \cref{eq:micro} numerically, the problem is discretized in
space using a second order accurate finite difference scheme, based
on a regular grid with mesh size $\Delta \xi = 2 \mu'/N_\mathrm{mic}$.
Let $a_{i, j} := a(-\mu' + i \Delta \xi,-\mu' + j \Delta \xi)$ and $\m_{i, j}(\tau) \approx \m^\varepsilon((-\mu' + i \Delta \xi,-\mu' + j \Delta \xi), \tau)$,
for $i, j = 0, ..., N_\mathrm{mic}$. Then the semi-discrete system obtained by discretization of \cref{eq:micro2} in space is
\begin{align*}
  \partial_\tau \m_{i,j}(\tau) = \f_{i,j}(\m), \quad \text{for } i, j = 1, ..., N_\mathrm{mic}-1,
\end{align*}
where
\begin{align*}
  \f_{i,j} (\m) :&= - \m_{i,j} \times \H_{i,j}(\m) - \alpha \m_{i,j} \times \left[\m_{i,j} \times \H_{i,j}(\m)\right], \\
  \H_{i,j}(\m) :&= \frac{1}{\Delta \xi^2} \left[a_{i+\tfrac{1}{2}, j} \m_{i+1, j} + a_{i-\tfrac{1}{2}, j} \m_{i-1, j} +
a_{i, j+\tfrac{1}{2}} \m_{i, j+1} + a_{i, j-\tfrac{1}{2}} \m_{i, j-1} \right.\\
 &\hspace{2cm}\left.- (a_{i+\tfrac{1}{2}, j} + a_{i-\tfrac{1}{2}, j} + a_{i, j+\tfrac{1}{2}} + a_{i, j-\tfrac{1}{2}})   \m_{i, j}
  \right].
\end{align*}
At $i = 0$, $j = 0$, $i = N_\mathrm{mic}$ or $j = N_\mathrm{mic}$,
we have
$\m_\mathrm{i, j} = \m_\mathrm{init}(-\mu' + i \Delta \xi, -\mu'+ j
\Delta \xi)$.  For time integration, we use the midpoint
extrapolation method (MPE) \cite{mpe}, a second order accurate
integrator that is norm preserving without any projections, which
makes it suitable for the non-normalized initial data. Let
$\m_{i,j}^k \approx \m_{i,j}(\tau_k)$, where
$\tau_k = k \Delta \tau$ for $k = 0, ..., M_\mathrm{mic}$ and $\Delta \tau =
\eta/M_\mathrm{mic}$. Then the time integration is described by
\[\m_{i,j}^{k+1} = \m_{i, j}^{k} - \Delta \tau \frac{\m^{k+1}_{i,j} + \m^k_{i,j}}{2}
  \times \h_{i,j}^{k+1/2},\]
where
\[\h_{i,j}(\m) := \H_{i,j}(\m) + \m_{i,j} \times \H_{i,j}(\m)\]
and $\h^{k+1/2}$ is a second order extrapolation approximating
$\h_{i,j}(\tfrac{1}{2}(\m^k + \m^{k+1}))$,
\begin{align*}
  \h_{i,j}^{k+1/2} := \frac{3}{2} \h_{i,j}(\m^k) - \frac{1}{2} \h_{i,j}(\m^{k-1}).
\end{align*}
Note that with this time stepping scheme, one only obtains stable
solutions given that
\begin{align}\label{eq:stab_cond}
  \Delta \tau \le C \Delta \xi ^2,
\end{align}
for some constant $C$ independent of $\Delta \xi$ but dependent on
the damping parameter $\alpha$. A detailed discussion of a micro
problem similar to the one discussed here is given in
\cite{paper3}. As suggested there, also in this paper we use
artificial damping in the micro problem and choose $\alpha \approx 1$,
which leads to an improved constant in the error estimate
in \Cref{thm:paper2} and thus convergence of the approximation errors
for shorter final times $\eta$. The averaging kernels $k$ and $k^0$
in \cref{eq:flux} are chosen such that $p_x = p_t = 3$ and
$q_x = q_t = q = 7$. As we moreover have $\mu \sim \varepsilon$ and
$\eta \sim \varepsilon^2$, we can simplify the error estimate for
the periodic case given in \cref{eq:err_terms}, which together with
the additional error terms due to the micro problem setup yields
\begin{align}\label{eq:error}
  |\F - \bnabla \m_\mathrm{init}(0,0) \A^H| \approx
  C \left[ \varepsilon + \left(\frac{\varepsilon}{\mu}\right)^{q + 2} + \left(\frac{\varepsilon^2}{\eta}\right)^{q + 1}\right] + E_{\mu'} + E_\mathrm{norm} + \mathcal{O}(\Delta \xi^2).
\end{align}
The last term here, $\mathcal{O}(\Delta \xi^2)$, is due to the
discretization error introduced when solving the micro problem numerically as described above.

While the error estimate in \Cref{thm:paper2} is only proved for
periodic material coefficients, we find that the upscaling errors
still behave according to \cref{eq:error} for somewhat more general
coefficients.  To demonstrate this, and to investigate how the
choices of $\mu, \eta$ and $\mu'$ influence the approximation error,
we consider a numerical example with the locally periodic material coefficient
\begin{align}
  a^\varepsilon(x) = 1.1 + \frac{1}{2} \left[\sin(2\pi x_1/\varepsilon) + \sin(2\pi x_2/\varepsilon)\right]\cos(2\pi(x_1+x_2)).
\end{align}
The corresponding homogenized matrix at the fixed location $x = (0,0)$
can be determined by freezing the slow variables and computing
$\A^H$ according to \cref{eq:AH}. This can be used to obtain a
reference solution for one micro problem, given macro data on one
triangle with barycenter in $(0,0)$.

\begin{figure}[h!]
  \centering
  \begin{subfigure}[b]{.31\textwidth}
    \includegraphics[width=\textwidth]{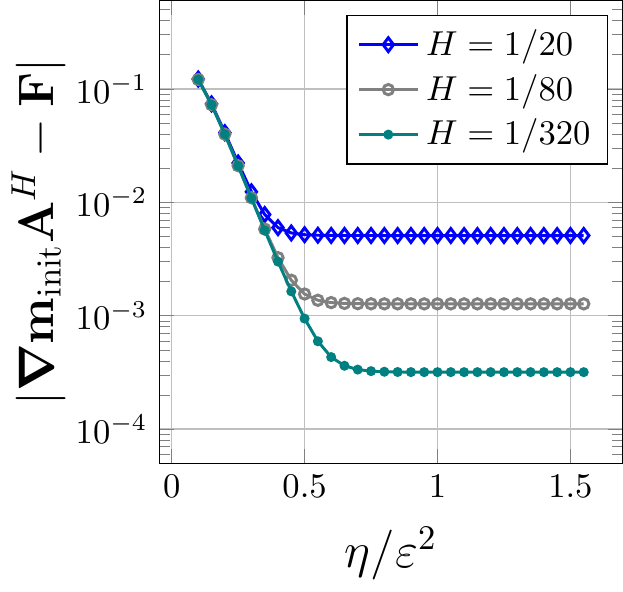}
    \caption{\centering Varying $\eta$, fixed  $\mu = 3.5\varepsilon$ and $\mu' = 15\varepsilon$. \label{fig:eta}}
  \end{subfigure}
  \quad
  \begin{subfigure}[b]{.31\textwidth}
    \includegraphics[width=\textwidth]{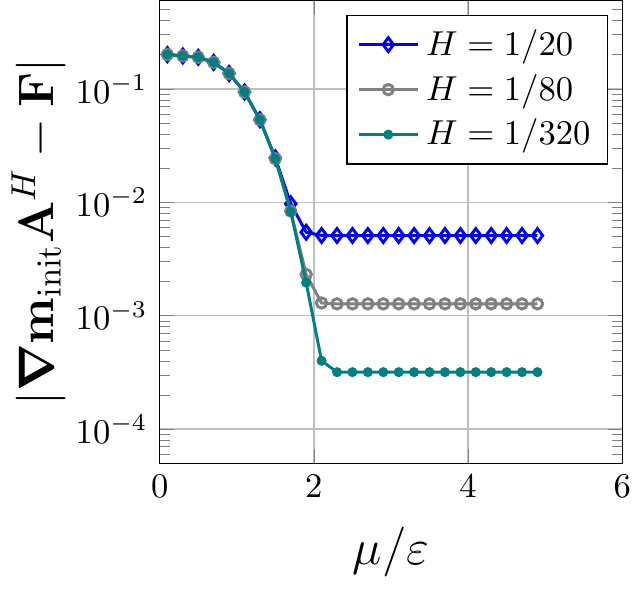}
    \caption{\centering Varying $\mu$, fixed $\eta = \varepsilon^2$ and $\mu'= 15\varepsilon$. \label{fig:mu}}
  \end{subfigure}
  \quad
  \begin{subfigure}[b]{.31\textwidth}
    \includegraphics[width=\textwidth]{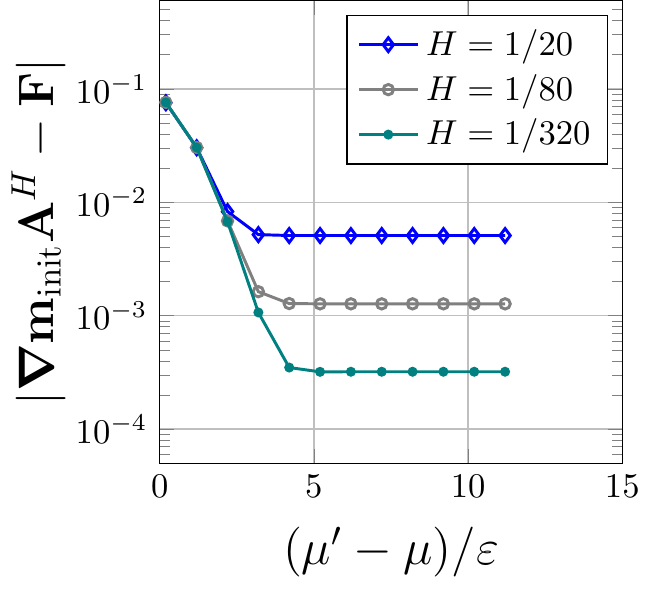}
    \caption{\centering Varying $\mu'$, fixed $\eta = \varepsilon^2$ and $\mu = 2.5\varepsilon$. \label{fig:box}}
  \end{subfigure}

  \caption{Influence of the micro problem parameters on the upscaling error.}
  \label{fig:micro_avg}
\end{figure}

To investigate the averaging errors, we then first fix
$\varepsilon = 10^{-4}$, and choose $\mu$ and $\mu'$ such that they
do not influence the upscaling error significantly, for this example
$\mu = 3 \varepsilon$ and $\mu' = 15 \varepsilon$. Moreover, we pick
$\Delta \xi$ such that the numerical discretization error is small
and select $\alpha = 1.5$, which means that we consider a setup with artificial
damping. Then the upscaling error is determined by $\eta$ and $H$,
and the estimate in \cref{eq:error} simplifies to
\begin{align*}
  |\F - \bnabla \m_\mathrm{init}(0,0) \A^H| \approx
  C \left(\frac{\varepsilon^2}{\eta}\right)^{q + 1} +  E_\mathrm{norm}.
\end{align*}
In \Cref{fig:eta}, the upscaling error for varying $\eta$ and
several choices of macro discretization length $H$ is shown. One can
observe that the error decreases rapidly as $\eta$ increases, until
it saturates at a level proportional to $H$, which corresponds to
$E_\mathrm{norm}$.

Next, we fix the length of the micro time
interval to $\eta = \varepsilon^2$ and vary the length of the
spatial averaging domain, $\mu$. All other parameters stay fixed.
In this case, the error estimate \cref{eq:error} becomes
\begin{align}
  |\F - \bnabla \m_\mathrm{init}(0,0) \A^H| \approx
  C \left(\frac{\varepsilon}{\mu}\right)^{q + 2} +  E_\mathrm{norm}.
\end{align}
The resulting errors are plotted in \Cref{fig:mu}. As shown there,
as $\mu$ becomes larger than $\varepsilon$, the error decreases rapidly and
then, once $\mu \gtrsim 2 \varepsilon$, saturates at a level
determined by the given value of $H$, the same as in \Cref{fig:eta},
corresponding to $E_\mathrm{norm}$.

Finally, both $\eta$ and $\mu$ are fixed and the micro domain size
$\mu'$ is varied. In this case, $E_{\mu'}$, for which no explicit
formula is known, and again $E_\mathrm{norm}$ determine the
resulting upscaling error, which is shown in \Cref{fig:box}. Based
on this plot, we come to the conclusion that taking
$\mu' \approx \mu + 5 \varepsilon$ is sufficient to reduce the error
introduced by the artificial boundary conditions so that it does not
influence the overall error in a significant way.

Further examples of how the upscaling error is influenced by the
choice of $\eta$, $\mu$ and $\mu'$ are given in \cite{paper3} for a
similar but slightly different micro problem.  We here conclude that
choosing $\eta$, $\mu$ and $\mu'$ large enough results in upscaling
errors determined only by $H$. Given a smaller value of $H$, lower
errors can be obtained by selecting the parameters $\eta$ and as a
consequence $\mu'$ larger. The optimal value for $\mu$ is only affected marginally.

The computational cost per micro problem can be described in terms
of the number of points used to resolve the $\varepsilon$ scale. Let
$P := \varepsilon / \Delta \xi$, then the total number of
discretization points for the micro problem is determined by
$N_\mathrm{mic}^2 = \left(2 \mu'/\Delta \xi\right)^2 = 4 P^2
\left(\mu'/\varepsilon \right)^2$. Moreover, due to
the time step restriction \cref{eq:stab_cond}, the number of time
steps per micro problem has to be chosen to be
$M_\mathrm{mic} = \eta/\Delta \tau = C \eta/\Delta \xi^2 = C P^2 \left(\eta /
  \varepsilon^2\right)$.  Hence, the overall cost per micro problem
is given by
\begin{align}
  \mathrm{cost~} \sim ~ C M_\mathrm{mic} N_\mathrm{mic}^2 \sim C P^4 \left(\frac{\mu'}{\varepsilon}\right)^2 \left(\frac{\eta}{\varepsilon^2}\right).
\end{align}
As $\mu' \sim \varepsilon$ and $\eta \sim \varepsilon^2$, this shows
that the computational cost per micro problem is independent of
$\varepsilon$. This makes it possible to apply the HMM approach even for very
small values of $\varepsilon$, where the resolution of the fine
scale with a conventional approach would result in tremendously high
computational cost.

Note, though, that to keep the overall cost down, one should choose
$\mu'$ and $\eta$ not larger than necessary and select a relatively
low value of $P$. For the example problems discussed in the next
section, $P \approx 10$ provides a reasonable compromise between
accuracy and computational cost.

\section{Numerical examples}\label{sec:num_ex}

In this section, we give several numerical examples to illustrate
the proposed scheme. We mostly focus on periodic material
coefficients to be able to provide corresponding homogenized solutions for reference. In the final
example, the locally periodic case discussed in \Cref{sec:upscaling}
is considered. In all the presented examples, we use homogeneous
Neumann boundary conditions as in \cref{eq:hom},
\cref{eq:macro_strong}.

\subsection{Circular domain example}\label{sec:ex1}
As a first numerical example, we consider a variation of the problem
suggested in \cite{alouges2008new}, with a circular 2D domain
$\Omega = B(0, 1)$ and initial data
\begin{align}\label{eq:ex_alouges_init}
  \M_\mathrm{init}(x) =
  \begin{bmatrix}
  -\frac{x_2}{r} \sin \left(\frac{\pi r}{2}\right), &
  \frac{x_1}{r} \sin \left(\frac{\pi r}{2}\right), &
  \cos\left(\frac{\pi r}{2}\right)
    \end{bmatrix}^T,
\end{align}
where $r = \sqrt{x_1^2 + x_2^2}$. Only the exchange term in the
effective field is considered,
$\H^\varepsilon = \bnabla \cdot (a^\varepsilon \bnabla \m^\varepsilon)$, where
$a^\varepsilon$ is assumed to be a periodic material coefficient,
\begin{align}\label{eq:mat_coeff1}
  a^\varepsilon(x) :=
  (1.1 + 0.25 \sin(2 \pi x_1/\varepsilon))(1.1 + 0.25 \sin(2 \pi x_2/\varepsilon)) + 0.7 \cos(2 \pi (x_1 - x_2)/\varepsilon).
\end{align}
The corresponding homogenized coefficient matrix, which is used to
solve the related homogenized problem for reference, is
\[\A^H \approx
  \begin{bmatrix}
    1.057 & 0.118 \\
    0.118 & 1.057
  \end{bmatrix},
\]
which is computed numerically with high accuracy for the simulation.
Moreover, the average of the given material coefficient,
\cref{eq:mat_coeff1}, is $a_\mathrm{avg} = 1.21$. We include the
solution to \cref{eq:LL} with the material coefficient replaced by
its average, $a_\mathrm{avg}$, in the example since this can be seen
as a naive approach to dealing with the oscillations in the
problem. We do not expect this to give the correct solutions. The
main reason of including it is to show that the chosen example is
relevant in the sense that its solution with the naive coefficient $a_\mathrm{avg}$
differs significantly from the correct solution.

\begin{figure}[!h]
  \centering
  \includegraphics[width=.9\textwidth]{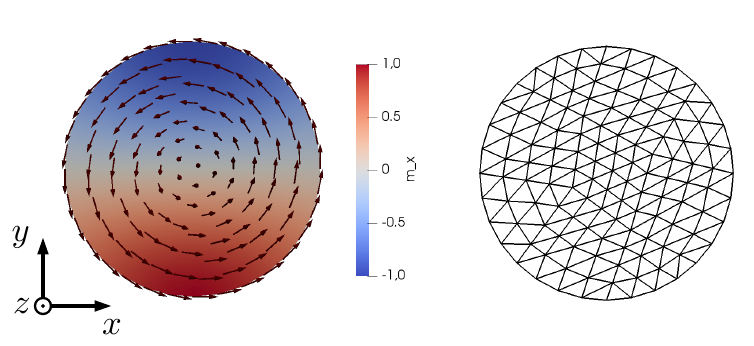}
  \caption{Initial data for the numerical example according to
    \cref{eq:ex_alouges_init} and computational mesh for the macro domain.
  }
  \label{fig:circle_init}
\end{figure}

The implementation of the finite element code for the numerical
examples in this article is done using the FEniCS project
\cite{logg2012automated}. For mesh generation, we use gmsh
\cite{geuzaine2009gmsh}. The initial data to the problem according
to \cref{eq:ex_alouges_init} and a mesh for the domain are shown in
\Cref{fig:circle_init}. There the domain is colored according to the
$x$-component of the data and the vectors show the direction of the
magnetization. Note that the $z$-direction is out of plane.

\begin{figure}[h!]
  {\hspace{1.8cm}\large \bf{Hom.} \hspace{2.3cm} \bf{HMM} \hspace{2.3cm} \bf{Avg.}}

  \begin{subfigure}[b]{\textwidth}
      \centering
    \includegraphics[width=\textwidth]{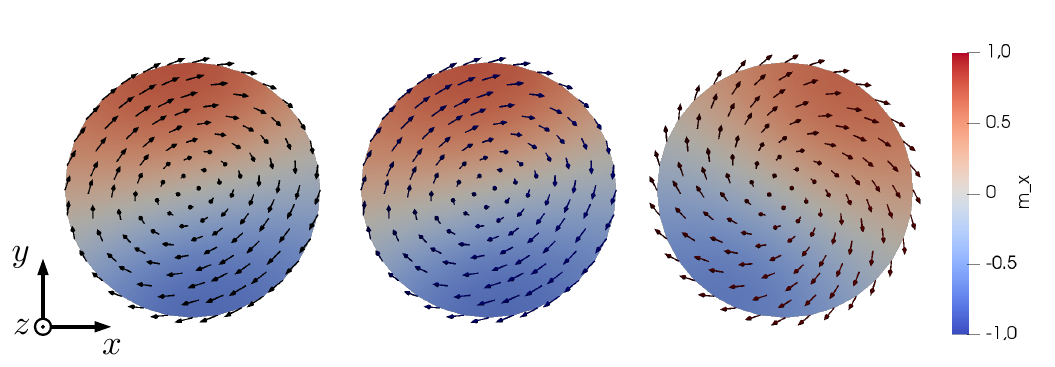}
  \end{subfigure}

  \caption{Numerical solution to homogenized problem (left), HMM
    solution (middle) and average coefficient solution (right) at
    time $T = 1.0$ to \cref{eq:LL} with material coefficient
    \cref{eq:mat_coeff1}, initial data according to
    \cref{eq:ex_alouges_init} and $\alpha = 0.2$ on the macro
    scale. Domain colored according to $x$-component of $\M$.
  }
  \label{fig:circle}
\end{figure}

For the HMM solution of this problem, we use artificial damping and
set $\alpha = 1.2$ in the micro problem. We choose
$\varepsilon = 10^{-4}$ in the example. Moreover, we consider two
different setups of averaging parameters:
\begin{itemize}[beginpenalty=10000,midpenalty=10000]
\item setup 1: $\mu' = 4.8 \varepsilon$, $\mu = 2.8 \varepsilon$ and
  $\eta = 0.45 \varepsilon^2$, which results in an upscaling error
  of approximately $0.002 + H$ for the micro problem around $(0,0)$,
\item setup 2: $\mu' = 3.25 \varepsilon$, $\mu = 2.1 \varepsilon$
  and $\eta = 0.15 \varepsilon^2$, which gives an upscaling error of
  approximately $0.02 + H$ for the same micro problem.
\end{itemize}
Consider a final time $T = 1.0$ and the macro discretization using
the mesh in \Cref{fig:circle_init}. Then in \Cref{fig:circle} the
corresponding solution to the homogenized problem, the HMM solution
with the micro parameters set according to setup 1, and the solution
when using the average of the material coefficient are shown. One
can observe that homogenized and HMM solution agree very well while
the averaged-coefficient solution is different.

Furthermore, in \Cref{tab:ex1}, the $L^2$ norm of the difference
between HMM solution and homogenized solution on a finer grid at
time $T = 0.5$ are given. With setup 1, corresponding to a
sufficiently low upscaling error, we observe second order
convergence of the HMM solution towards the reference solution. With
setup $2$, on the other hand, the errors initially decrease when
refining the computational grid but then seem to saturate. This is
due to the fact that in this case, the upscaling error is larger and
hence affects the HMM solution more. We expect a similar effect to
appear also with the first micro parameter setup, but for more refined grids, with lower
$H_\mathrm{min}$.

\renewcommand{\arraystretch}{1.15}
\vspace{.1cm}
\begin{table}[h!]
  \centering
  \begin{tabular}[h]{l || c | c|| c| c|  }
    $H_\mathrm{min}$ & $\|\M - \M_\mathrm{ref}\|_{L^2}$ & conv. & $\|\M - \M_\mathrm{ref}\|_{L^2}$ & conv. \\
    & setup 1 & order & setup 2 & order\\
        \hline
    0.684            & 0.0663 & & 0.0858&  \\
    0.342            & 0.0197 & 1.75 & 0.0273 & 1.65\\
    0.171            & 0.00557 & 1.82 & 0.0129 & 1.08 \\
    0.0855           & 0.00139 & 2.00 &0.0101 & 0.353\\
  \end{tabular}
  \caption{Convergence of the HMM solution at time $T=0.5$ with
    respect to a homogenized reference solution on a fine grid. Same
    choice of initial data and material coefficient as in \Cref{fig:circle}
    and again $\alpha = 0.2$.}
  \label{tab:ex1}
\end{table}

\begin{figure}[!htb]
  \centering
  \includegraphics[width=\textwidth]{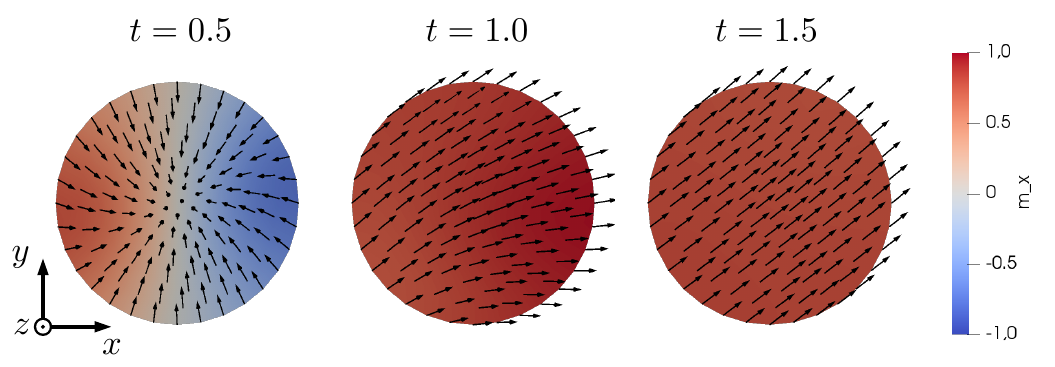}
  \caption{Influence of an external field. Same example as in \Cref{fig:circle} but with $\H_\mathrm{ex} = [10, 10, 0]^T$ applied for $t > 0.5$.}
  \label{fig:ex1_field}
\end{figure}

We then furthermore extend the example problem and additionally consider an applied external
field,
\[\H_\mathrm{ex} = \left[10, 10, 0\right]^T,\]
which is turned on at time $T = 0.5$. This causes the magnetization
vectors to align accordingly, as shown in \Cref{fig:ex1_field}.

\clearpage

\subsection{Ring example}
\begin{figure}[hp!]
  \centering
  \includegraphics[width=.9\textwidth]{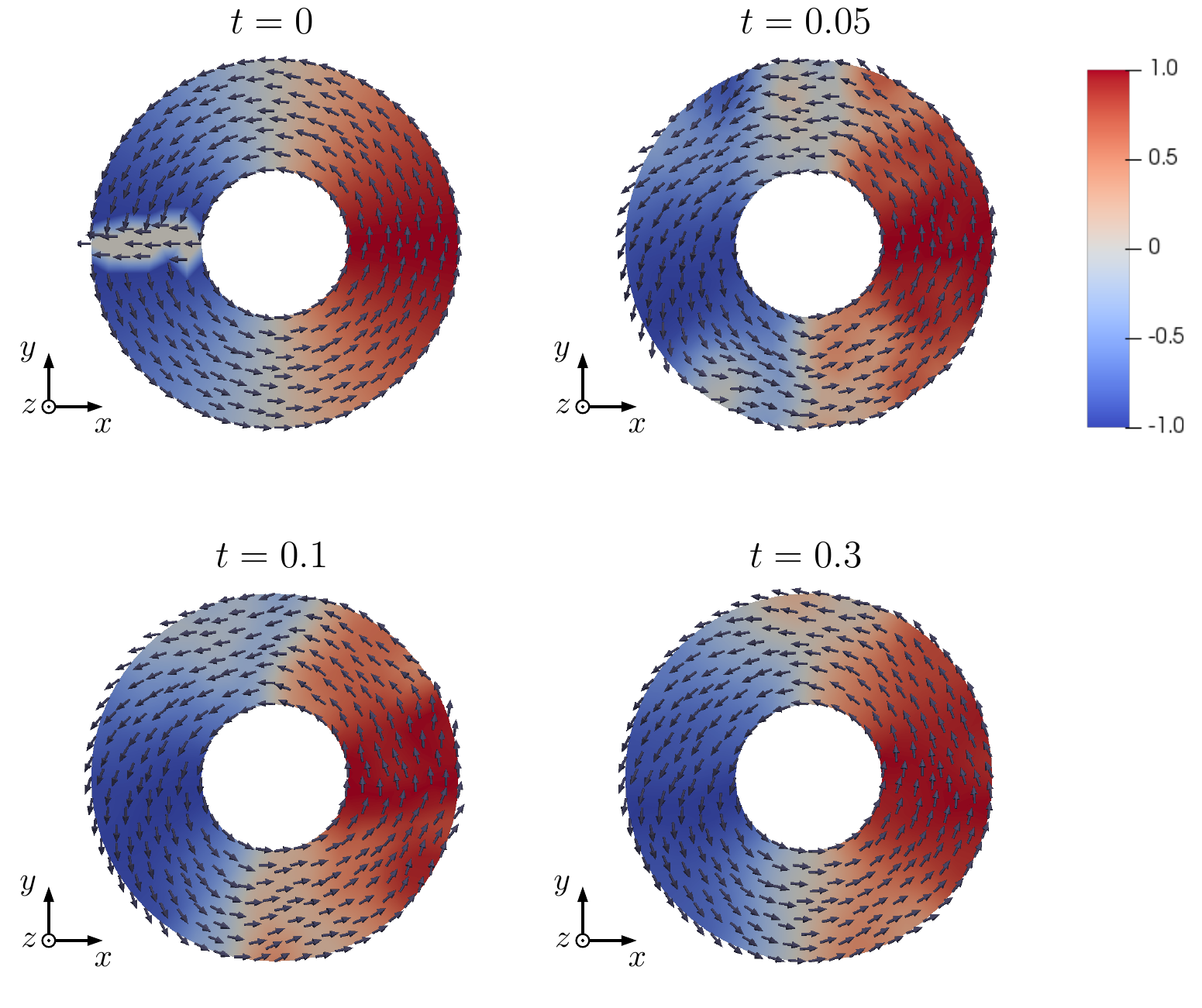}
  \caption{HMM solution to the ring example problem. Domain colored according to $y$-component of $\M$.}
  \label{fig:ring}
\end{figure}

Inspired by \cite{chaves2009micromagnetic}, we next consider a
ring-shaped domain with outer radius $R_\mathrm{out} = 1$ and inner radius
$R_\mathrm{in} = 0.4$. Initially, the magnetization in the ring is given by
\begin{align}\label{eq:init_ring}
  \M_\mathrm{init}(x) =
  \begin{bmatrix}
  -\frac{x_2}{r} , &   \frac{x_1}{r} , & 0
    \end{bmatrix}^T,
\end{align}
where again $r = \sqrt{x_1^2 + x_2^2}$,
which is close to a stable state, except for a small subsection of
the domain where we instead set $\M_\mathrm{init} = [1, 0, 0]^T$, a
simplified representation of a defect. This is shown in
\Cref{fig:ring} at time $t = 0$. For this example problem, we set
the damping parameter to $\alpha = 0.02$ (on the macro scale) and
use the same material coefficient $a^\varepsilon$ as in the previous
example, \cref{eq:mat_coeff1}. Only the exchange-term is considered
in the effective field.  The micro parameters are selected to be
$\eta = 0.3 \varepsilon^2$, $\mu = 2.1 \varepsilon$ and
$\mu' = 4.25$, in between the two setups considered in the previous subsection, and
$\varepsilon = 10^{-3}$.  Furthermore, artificial damping with
$\alpha = 1.2$ is used in the micro problem.

The resulting HMM solution at several points in time is shown in
\Cref{fig:ring}. Due to the defect, the magnetization throughout the
ring is affected and varies as the magnetization vectors strive
towards a state of alignment.

\begin{figure}[hp!]
  \centering
  \includegraphics[width=.87\textwidth]{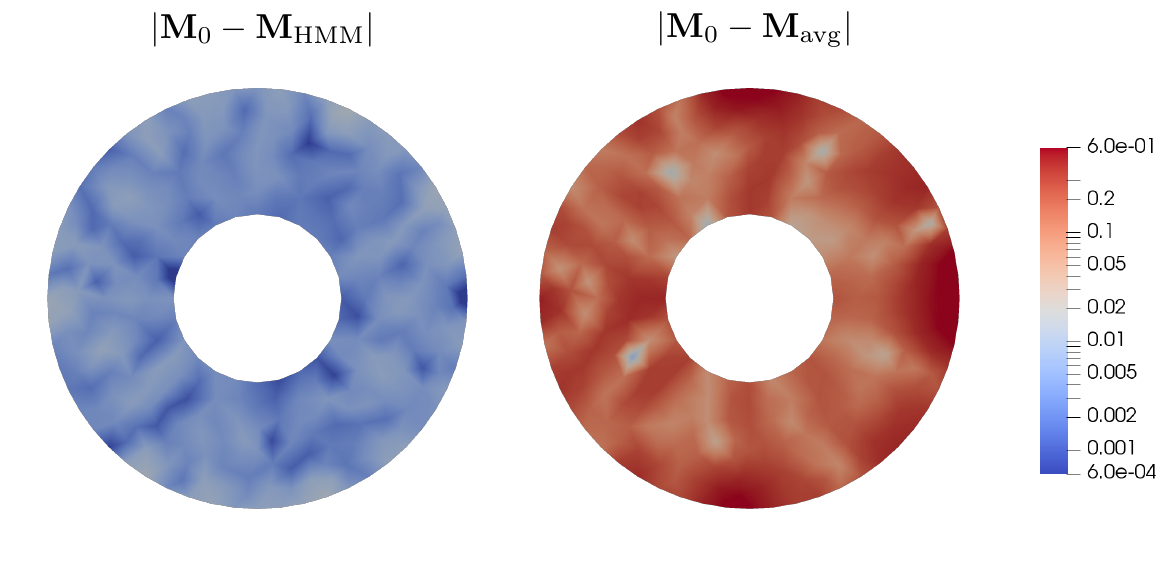}
  \caption{Difference between homogenized solution (for reference) and HMM solution or average coefficient solution, respectively.}
  \label{fig:ring_err}
  \vspace{1cm}
\end{figure}

Furthermore, for the final time $T = 0.3$, the error between
homogenized solution $\M_0$ and HMM solution as well as between
homogenized solution and the solution obtained when using the
average of the material coefficient is shown in \Cref{fig:ring_err}.
One can observe that despite the coarser choice of micro problem
parameters, the HMM solution agrees well with the homogenized
solution, while the average coefficient solution shows high errors.

\subsection{Landau-state example}
As a further example, we consider relaxation towards a so-called
Landau-state.
A domain of length and width 250 nm and thickness 2 nm is given,
which in a non-dimensional setting is represented by a 2D
unit-square domain. The third dimension, with much smaller extent
than the other two, is only considered in the computation of the
demagnetization term.  A time interval of 1 ns is assumed and the
damping parameter is selected as $\alpha = 0.05$.
The initial data is set such that the magnetization
eventually attains a so-called Landau state,
\[\M_\mathrm{init} =
  \begin{cases}
    [0, ~~1, ~0]^T & x_1 < 0.5 \\
    [0, -1, ~0]^T & \mathrm{otherwise}.
  \end{cases}
\]
The effective field in
this example consists of exchange term and demagnetization,
\[\H^\varepsilon = \bnabla \cdot (a^\varepsilon \bnabla \m^\varepsilon) + \H_\mathrm{dem}.\]
In general, the demagnetization field is given by
$\H_\mathrm{dem} = - \nabla u$, where the scalar potential $u$ solves
\[ \Delta u =
  \begin{cases}
    \bnabla \cdot \M & \text{in~} \Omega, \\
    0 & \text{in~} \Real^3 \backslash \Omega.
  \end{cases}
\]
  For reasons of simplicity, the demagnetization term for this
  example problem is computed using the algorithm proposed in
  \cite{newell1993generalization}, specifically using a C++
  implementation along the lines of \cite{70lines}. This algorithm
  uses the fact that $\H_\mathrm{dem}$ can also be expressed as
  \[\H_\mathrm{dem} = \int_\Omega \tilde \N(x - y) \M(y) dy, \quad \text{~where~} \quad \tilde \N(x -y) = -\frac{1}{4 \pi} \bnabla \bnabla_y \frac{1}{|x-y|}.\]
  For an efficient implementation, the demagnetization tensor
  $\tilde \N$ is pre-computed and the convolution integral is
  evaluated using the Fast Fourier transform.  It is computed on a
  regular grid in a separate module which is coupled with the finite
  element solver.

\begin{figure}[h!]
  \centering
  \includegraphics[width=\textwidth]{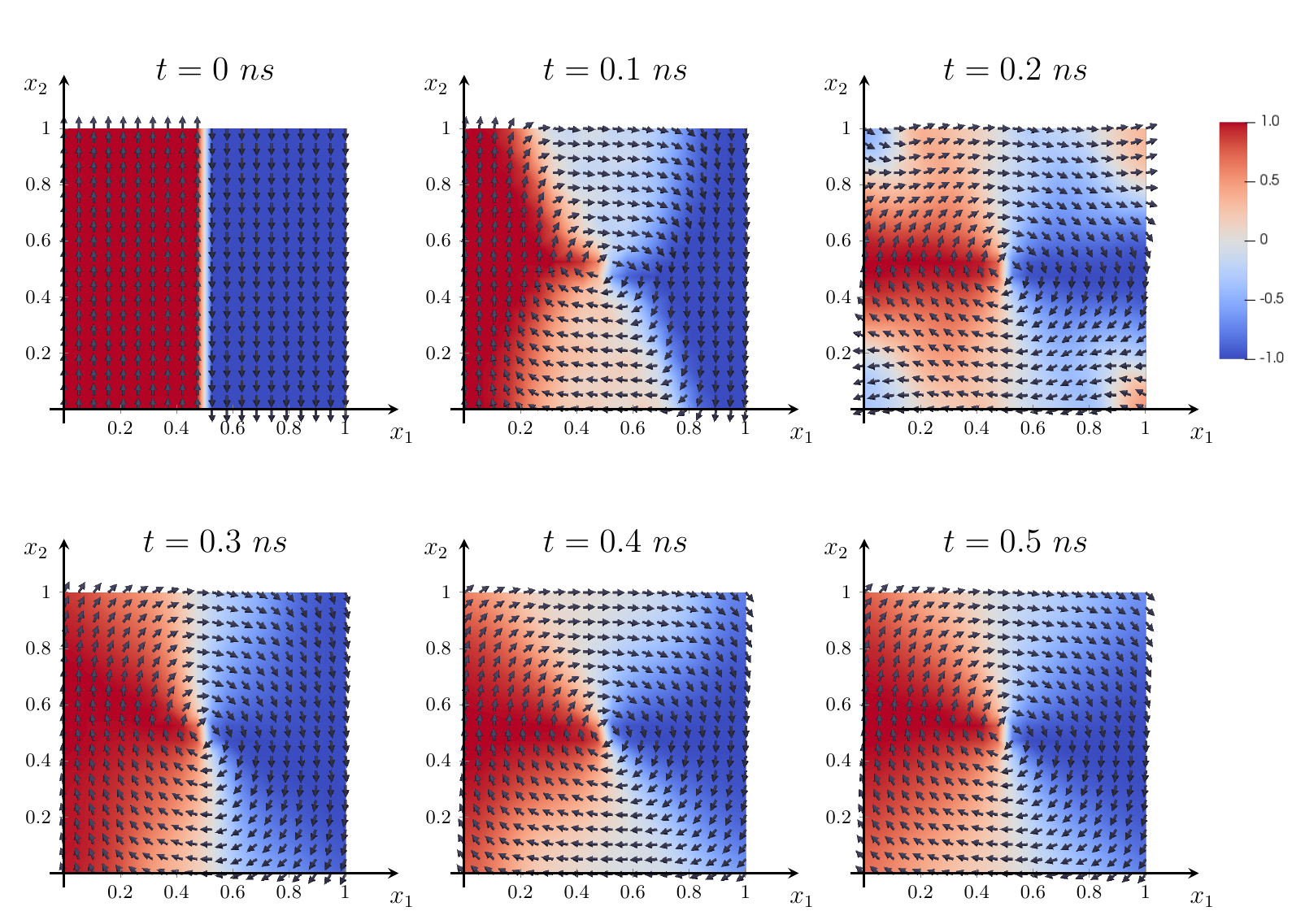}
  \caption{HMM solution to the Landau-state example problem for different points in time. Domain colored
    according to $y$-component of the magnetization.}
  \label{fig:landau}
\end{figure}

To make it possible to compare to a homogenized solution for
reference, we use another periodic material coefficient,
\[a^\varepsilon = \exp\left[\cos(2 \pi (x+y)/\varepsilon) - 0.25 \sin(2 \pi x/\varepsilon)\right].\]
The corresponding homogenized matrix is approximately
\[
  \A^H =
  \begin{bmatrix}
     1.014  &-0.234\\  -0.234 &   1.04
  \end{bmatrix}.\]
In this example problem, we set $\varepsilon = 10^{-3}$ and select
the micro problem parameters to be $\eta = 0.3 \varepsilon^2$, $\mu = 2.1 \varepsilon$
and $\mu' = 5.5 \varepsilon$. For the micro problem, $\alpha = 1.2$. The obtained HMM solution
at times corresponding to $0$ to $ 0.5$ ns is shown in \Cref{fig:landau}.

\begin{figure}[h]
  \centering
  \includegraphics[width=.8\textwidth]{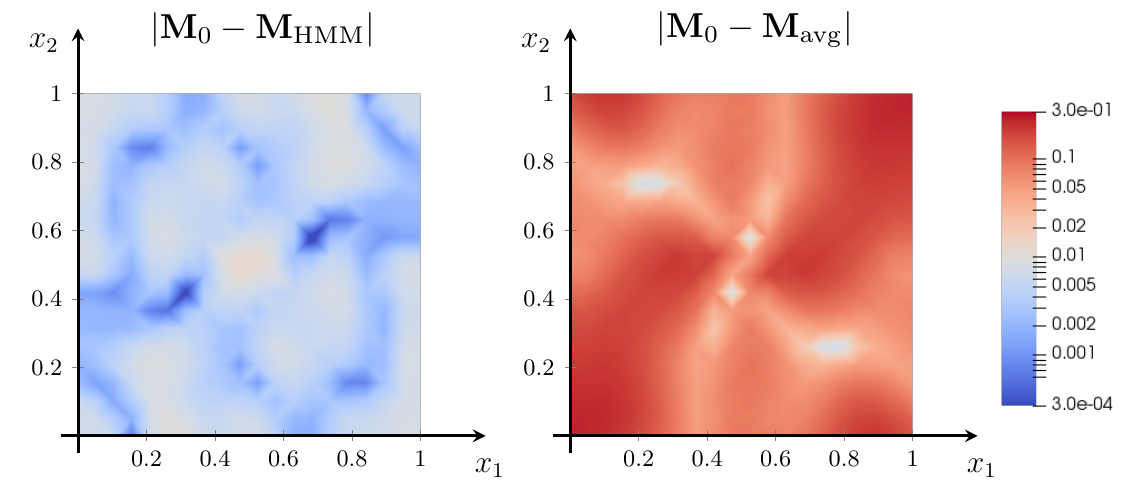}
  \caption{Error in HMM and avg. coefficient solution at time $t = 0.5 ns$ }
  \label{fig:err_landau}
\end{figure}

Moreover, the difference between a corresponding homogenized
solution $\M_0$ and the HMM solution at time corresponding to
$t = 0.5$ ns is shown in \Cref{fig:err_landau}. For comparison, also
the error between the homogenized solution $\M_0$ and what is
obtained when approximating $a^\varepsilon$ by its average
$a_\mathrm{avg}$, a solution $\M_\mathrm{avg}$, is shown. One can
observe that the HMM solution agrees considerably better with the
homogenized solution than $\M_\mathrm{avg}$.

\subsection{Modification of $\mu$MAG standard problem 4}
\begin{figure}[h!]
  \centering
  \includegraphics[width=\textwidth]{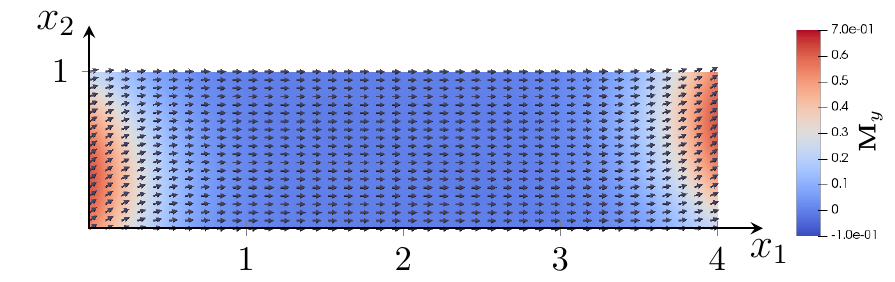}\qquad\qquad
  \caption{The equilibrium s-state obtained from the initial data
    \cref{eq:init_prob4} with material coefficient $a^\varepsilon$
    given by \cref{eq:aex_loc}, colored according to $y$-component
    of $\M$.}
  \label{fig:s-state}
\end{figure}

As a final example, we consider a variation of the $\mu$MAG
standard problem 4 \cite{mumag}, with a locally periodic material coefficient,
\begin{align}\label{eq:aex_loc}
  a^\varepsilon(x) = 1.1 + \frac{1}{2}\left[\sin(2\pi x_1/\varepsilon) + \sin(2\pi x_2/\varepsilon)\right]\cos(2\pi(x_1+x_2)),
\end{align}
the same coefficient as in the example in \Cref{sec:upscaling}. In the $\mu$MAG problem, a thin film domain
of length 500 nm, width 125 nm and thickness 3 nm is considered. The
dimensional form of the Landau-Lifshitz equation is used, with
material parameters similar to permalloy.
In this article, we use a corresponding non-dimensionalized,
rescaled setup and a 2D domain of length 4 and width 1. The third
dimension is again only
taken into account in the computation of the demagnetization. As initial
data for this problem, a so-called equilibrium s-state is used, as
shown in \Cref{fig:s-state}. We obtain this s-state by first setting
the magnetization to
\begin{equation}
  \label{eq:init_prob4}
  \M_\mathrm{init} =
  [\cos(\pi x_1/8), ~ \sin(\pi x_1/8), ~ 0]^T,
\end{equation}
 and then relaxing until equilibrium is reached.
Once the initial s-state is achieved, an external field is applied,
corresponding to $[-24.6, 4.3, 0.0]^T/\mu_0 \mathrm{\,mT}$, where
$\mu_0$ denotes the vacuum permeability,
$\mu_0 = 4\pi \cdot 10^{-7} \mathrm{\,N/A^2}$, in the dimensional
problem. The damping parameter is set to $\alpha = 0.02$.

In this final example problem, we consider exchange interaction, applied
external field and demagnetization,
\[\H^\varepsilon = \bnabla \cdot (a^\varepsilon \bnabla \m^\varepsilon) + \H_\mathrm{ex} + \H_\mathrm{dem}.\]
The demagnetization is computed as described for the previous example.

  The HMM micro problem parameters for this example are chosen
  according to the discussions in \Cref{sec:upscaling}, based on
  \Cref{fig:micro_avg}. We use $\alpha = 1.5$ and set
  $\eta = 0.5 \varepsilon^2$, $\mu = 2.1 \varepsilon$ and
  $\mu' = 6.5 \varepsilon$, where in this example, $\varepsilon = 10^{-3}$.

\begin{figure}[h!]
  \centering
  \includegraphics[width=\textwidth]{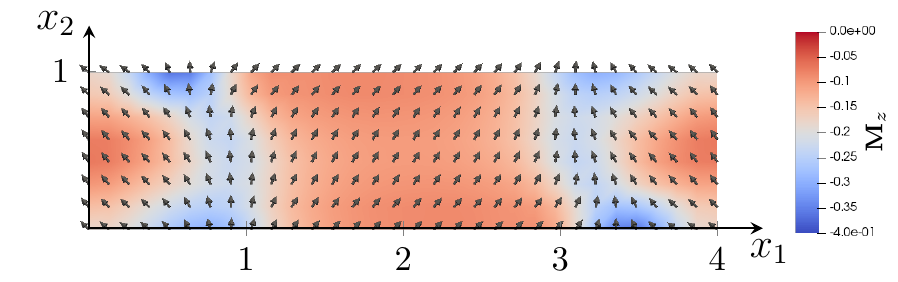}
  \caption{HMM solution when the average of the $x$-component,
    $\langle M_x \rangle$, first crosses zero. Domain colored
    according to $z$-component of the magnetization.}
  \label{fig:x_zero}
\end{figure}

  The corresponding HMM solution at the time when the average of the
  $x$-component of $\M$ first crosses zero is shown in
  \Cref{fig:x_zero}. It matches well with the general expected
  behavior based on the results reported for example at \cite{mumag}.

\section{Conclusion and outlook}

In this article, we presented a Heterogeneous Multiscale Method for
the Landau-Lifshitz equation with an oscillatory material coefficient
which combines a finite element macro model with a finite difference
micro problem based on results from \cite{paper2}. Numerical
examples show the viability of the approach. It is thus possible to
treat scale-separated problems with arbitrarily small variations
$\varepsilon$, given a certain structure of those variations, such
as local periodicity.

The current model includes the contributions of exchange interaction,
applied external field and demagnetization to the effective field.
However, in the given implementation, the demagnetization can only be
computed for rectangular domains. This should be adjusted by using a
different approach for solving the demagnetization problem in the
future.  Moreover, it would be interesting to study how to also
include anisotropy effects into the model.

To keep focus on the multiscale aspect of the problem and for
reasons of simplicity, we in this article chose a weak formulation
for the macro model that does not specifically comply with the
orthogonality structure given for the continuous problem. It should be rather
simple to combine the proposed approach with schemes designed to
take this into account such as explicit-in-time versions of
\cite{alouges2006fe,alouges2008new,alouges2014precise}. Possible effects of doing so could be studied in the future.

A further interesting direction for future research could be to
investigate whether it is possible to extend basis-representation
techniques that have been proposed for HMM for linear problems, for example
in \cite{engquist2009multi}, to the given Landau-Lifshitz
problem. This could drastically reduce the number of micro problems
that have to be solved and hence computational cost.

\subsection*{Acknowledgments} The authors would like to thank
Prof. Olof Runborg and Prof. Gunilla Kreiss for many useful
discussions along this work.

\bibliographystyle{acm}
\bibliography{hmm}

\end{document}